\documentclass{article}
\usepackage{amssymb,latexsym}
\usepackage{amsmath}
\usepackage{graphics}
\newtheorem{theorem}{Theorem}
\newtheorem{conjecture}{Conjecture}

\def \fo {f_{\text {opt}}}
\def \Pr {P_{3t+r}} \def \P0 {P_{3t}} \def \P1 {P_{3t+1}} \def \P2 {P_{3t+2}}
\def \ctr {C_{3t+r}}  \def \co{C_{3t+1}}  \def \cm {C_{3t-1}}

\begin{document}
\title{Optimal pebbling of paths and cycles}
\author{C. Wyels\\Dept. of Mathematics \& Physics\\
California Lutheran University\\Thousand Oaks, California 91360 \\
{\small\texttt{wyels@clunet.edu}}\\
\and 
T.~Friedman\\Dept. of Computer Science, Mathematics \& Statistics\\
Mesa State College\\
Grand Junction, CO  81501 \\
{\small\texttt{tfriedma@mesastate.edu}}
}
\date{May 30, 2003}

\maketitle

\begin{abstract}
Distributions of pebbles to the vertices of a graph are said to be solvable when a pebble may be moved to any specified vertex using a sequence of admissible pebbling rules. The optimal pebbling number is the least number of pebbles needed to create a solvable distribution.  We provide a simpler proof verifying Pachter, Snevily and Voxman's determination of the optimal pebbling number of paths, and then adapt the ideas in this proof to establish the optimal pebbling number of cycles.  Finally, we prove the optimal-pebbling version of Graham's conjecture.

\vspace{.05in}
\end{abstract}

\section{Introduction}\label{S:intro}
The concept of pebbling graphs was introduced by Lagarias and Saks to rephrase a number theoretic conjecture posed by Erd\"{o}s and Lemke \cite{Chung}.  The idea consists of distributing pebbles to the vertices of a graph, stating a rule for moving pebbles, and asking when at least one pebble may be moved to any vertex.  The rule governing moving pebbles states that one may remove two pebbles from a vertex and subsequently place one on an adjacent vertex.  We say that a distribution of pebbles to the vertices of a graph is solvable when each vertex can receive a pebble via a sequence of pebbling moves.  (This includes sequences of length 0, i.e.~when the initial distribution places at least one pebble on the vertex in question.)  Two types of pebbling numbers follow immediately:
\begin{enumerate}
\item The pebbling number of a graph $G$, $f(G)$, is the least number such that every distribution of $f(G)$ pebbles to the vertices of $G$ is solvable.
\item The optimal pebbling number of a graph $G$, $\fo(G)$, is the least number such that there exists a solvable distribution of $\fo(G)$ pebbles on $G$.
\end{enumerate}  

Chung established the first results in this area, giving some bounds for pebbling numbers as well as the pebbling numbers for hypercubes, paths, and complete graphs \cite{Chung}.  Others have established the pebbling numbers for various classes of graphs, e.g. cycles \cite{SnevilyFoster,Herscovici}, complete bipartite graphs \cite{XL}, and trees \cite{Moews92}.  Pebbling numbers of products of graphs have also been studied:  we discuss this briefly in Section 4.  Optimal pebbling numbers are not as well documented:  formulas are known only for paths \cite{PSV}, complete $m$-ary trees \cite{FuShiue}, hypercubes \cite{Moews98}, and graphs of diameter 2.  Interestingly, exact values for optimal pebbling numbers are known only for paths, complete $m$-ary trees with $m \geq 3$, and diameter-2 graphs.  In Section 3 we establish an exact value for the optimal pebbling number of cycles.

To enhance the readability of our proofs, we offer the following notational and linguistic conveniences.
\begin{itemize}
\item We often refer to a distribution by name.  For the distribution $D$ we use $|D|$ to indicate the size of $D$, i.e.~the number of pebbles in $D$.  
\item Suppose the vertices $v_i,v_{i+1},\dots ,v_{i+j}$ form a path within a graph $G$.  We write $D(\lbrack v_i,v_{i+1},\dots ,v_{i+j} \rbrack ) = \lbrack p_i,p_{i+1},\dots p_{i+j} \rbrack$ to indicate that $D$ places $p_k$ pebbles on vertex $v_k$ for each $k = i, i+1, \dots ,i+j$.  Similarly, $D(v) = p$ means that $D$ places $p$ pebbles on the vertex $v$.  We then say that a vertex $v$ is {\it occupied} if $D(v)>0$ and {\it unoccupied} when $D(v)=0$.
\item We often alter a graph by {\it removing a vertex} of degree 1 or 2.  If we remove a vertex $v$ with $\deg(v)=1$ then it is understood that the edge incident to $v$ is also removed.  If $\deg(v)=2$ then the two edges incident to $v$ are replaced with an edge incident to both of $v$'s original neighbors.  Thus a path (respectively, a cycle) on $n$ vertices from which a vertex is removed becomes simply a path (cycle) on $n-1$ vertices.  
\end{itemize}

\section {Optimal Pebbling of Paths}
Pachter, Snevily and Voxman establish the optimal pebbling number for path graphs \cite{PSV}.  The formula for $\fo(P_n)$ depends on the value of $n \mod 3$, so we write $P_n$, the path on $n$ vertices, as $P_n=\Pr$, where $r \in \{0,1,2\}$.
\begin{theorem} \cite{PSV} The optimal pebbling number of the path on $3t+r$ vertices is $2t+r$, i.e. $\fo(\Pr)=2t+r$ .
\end{theorem}
{\bf Proof.}
We rephrase the distribution given in \cite{PSV} to show that $\fo(\Pr)\leq 2t+r$:  Label the vertices of $\Pr$ sequentially as $v_1, v_2, \dots v_{3t+r}$.  For $1 \leq i \leq 3t$, set $D(v_i) = 2$ if $i \equiv 2 \mod 3$, and $D(v_i)=0$ otherwise.  If $r\geq 1$ set $D(v_{3t+1})=1$, and if $r=2$ then set $D(v_{3t+2})=1$ as well.  Observe that each vertex is either occupied or adjacent to a vertex with two pebbles.  Thus the distribution is solvable.

We ask the reader to check by hand that $\fo(\Pr) \geq 2t+r$ when $3t+r$ is 1, 2, and 3.  To show $\fo(\Pr) \geq 2t+r$ we assume the contrary, and that the smallest index for which the theorem fails is $3t+r$.  

\noindent {\bf Case 1:  $r = 0$.}
\vspace{.03in}
\newline
(The logic used here is that of \cite{PSV}.)  We rewrite $P_{3t-1}$ as $P_{3(t-1)+2}$.  Since this path has a smaller index than the least index for which the theorem fails, we are assured that $\fo(P_{3(t-1)+2}) = 2(t-1)+2 = 2t$.  From our assumption we know that $\fo(P_{3t}) < 2t$.  As $\fo(P_{3t-1}) \leq \fo(P_{3t})$, we have 
$$2t = \fo(P_{3t-1}) \leq \fo(P_{3t}) < 2t,$$
an impossibility.

\noindent {\bf Case 2:  $r = 1$ or $r=2$.}
\vspace{.03in}
\newline  Since we are assuming that $\fo(\Pr)\leq 2t+r-1$, we may choose a solvable distribution $D$ of $\Pr$ such that $|D|=2t+r-1$.  We will modify $D$ to create $D^*$, a solvable distribution of $P_{3t+r-1}$ with fewer than $\fo(P_{3t+r-1})$ pebbles, thus producing our desired contradiction.

{\it Subcase 2.1}:  $\Pr$ contains a vertex $v$ with $D(v)=1$.
\vspace{.03in}
\newline
We modify $D$ and $\Pr$ by removing $v$ and the pebble on it, thus creating a distribution $D^*$ of $2t+r-2$ pebbles on $P_{3t+r-1}$.  We claim that all the remaining vertices may still be pebbled, in other words, that $D^*$ is solvable.  Suppose $v$ was involved in a pebbling step needed to establish the solvability of $D$.  Then $v$ is not an endpoint of the path.  Number the vertices of $\Pr$ sequentially and let $v=v_i$.  Without loss of generality we may assume that vertices indexed by $j$, $j>i$, required pebbles moved from $v_{i-1}$ through $v_i$ to $v_{i+1}$ to pebble them.  Assume that, using distribution $D$, $a$ pebbles could be moved from $v_{i-1}$ to $v_i$.  Then $\lfloor \frac{a+1}2  \rfloor$ pebbles could be moved from $v_{i-1}$ through $v_i$ to $v_{i+1}$.  Using $D^*$, the same $a$ pebbles were moved from $v_{i-1}$ to $v_i$ may now be moved directly from $v_{i-1}$ to $v_{i+1}$.  As $a \geq \lfloor \frac{a+1}2  \rfloor$ for all $a \geq 1$, we see that all vertices initially reachable from $D$ are still reachable from $D^*$.  Of course, any vertex that was reachable from $D$ without using $v$ initially is still reachable from $D^*$.  

{\it Subcase 2.2}: $D$ puts at least 2 pebbles on all occupied vertices of $\Pr$.
\vspace{.03in}
\newline
We label the vertices of $\Pr$ sequentially as $v_1,v_2,\dots ,v_{3t+r}$.  Let $i$ be the smallest index for which $v_i$ is occupied and $v_{i+1}$ is unoccupied.  (If no such $i$ exists then renumber the vertices of $\Pr$ starting from the opposite end.)  Create $D^*$ by removing vertex $v_{i+1}$ from the graph and removing two pebbles from vertex $v_i$.  Additionally, if $i \neq 1$, place one additional pebble on $v_{i-1}$.  We now have $D^*$, a distribution with no more than $2t+r-2$ pebbles on $P_{3t+r-1}$.

As the two pebbles removed from $v_i$ could contribute exactly 1 pebble to pebbling steps used to pebble any vertex $v_h$, $h < i$, the additional pebble placed on $v_{i-1}$ ensures that any vertex $v_h$, $h < i$, that could be pebbled from $D$ may still be pebbled from $D^*$.  

It remains only to show that any vertices $v_{i+j}$ ($j \geq 2$) that could be pebbled from $D$ may still be pebbled from $D^*$.  Suppose $D(v_i)=a$;  then $D^*(v_i)=a-2$.  
Suppose that, starting from $D$, $b$ pebbles were collected on vertex $v_{i-1}$.  Then $\lfloor \frac b2 \rfloor$ pebbles could be moved to $v_i$, permitting $\lfloor \frac 12 (\lfloor \frac b2 \rfloor + a) \rfloor$ pebbles to be moved to $v_{i+1}$, and finally $\lfloor \frac 12 \lfloor \frac 12 (\lfloor \frac b2 \rfloor + a) \rfloor\rfloor$ pebbles to be moved to $v_{i+2}$.  
Now, starting from $D^*$, $\lfloor \frac {b+1}2 \rfloor$ pebbles may be moved from $v_{i-1}$ to $v_i$, and these $\lfloor \frac {b+1}2 \rfloor$ pebbles added to $v_i$'s $a-2$ pebbles allowing the placement of $\lfloor \frac 12 (\lfloor \frac {b+1}2 \rfloor + a-2) \rfloor$ pebbles on $v_{i+2}$.  It is straightforward to verify that $\lfloor \frac 12 (\lfloor \frac {b+1}2 \rfloor + a-2) \rfloor \geq \lfloor \frac 12 \lfloor \frac 12 (\lfloor \frac b2 \rfloor + a) \rfloor\rfloor$ whenever $a \geq 2$, so once again, all $v_{i+j}$ that could be pebbled starting from $D$ can still be pebbled from $D^*$.

We conclude that our new distribution $D^*$ is a solvable distribution on $P_{3t+r-1}$.  However, $|D^*|\leq 2t+r-2$, thus contradicting the fact that $\fo(P_{3t+r-1}) \linebreak= 2t+r-1$.  This discrepancy forces us to conclude that our initial assumption is false, and so the theorem holds for all positive integers.


\section{Optimal Pebbling of Cycles}
The proof given for path graphs extends nicely to provide the optimal pebbling number for cycle graphs, $C_n$.  Again, we are interested in $n\mod 3$, so we write $C_n$, the cycle on $n$ vertices, as $C_n=\ctr$, where $r \in \{0,1,2\}$.
\begin{theorem} The optimal pebbling number of the cycle on $3t+r$ vertices is $2t+r$, i.e. $\fo(\ctr)=2t+r$.
\end{theorem}
{\bf Proof.}  We begin by constructing a solvable distribution of size $2t+r$.  Number the vertices of $\ctr$ sequentially, $\ctr=(v_1, v_2, \dots, v_{3t+r})$.  For $1 \leq i \leq 3t$, set $D(v_i) = 2$ if $i \equiv 2 \mod 3$, and $D(v_i)=0$ otherwise.  If $r\geq 1$ set $D(v_{3t+1})=1$, and if $r=2$ then set $D(v_{3t+2})=1$ as well.    Observe that in each case, every vertex is either occupied or is adjacent to a vertex with two pebbles.  Thus, the distribution is solvable. 

It is straightforward to show that the theorem holds for cycles of lengths 3, 4, and 5.  It remains to show that $\fo(\ctr)\geq2t+r$ for $3t+r\geq 6$.  Suppose, to the contrary, that $3t+r$ is the least integer for which $\fo(\ctr)<2t+r$.

\noindent {\bf Case 1:  $r = 0$.}
\vspace{.03in}
\newline  Since $\cm=C_{3(t-1)+2}$, we have $\fo(\cm)=2(t-1)+2=2t$.  However, $\fo(\cm)\leq\fo(C_{3t})$.  Thus, $$2t=\fo(\cm)\leq\fo(C_{3t})<2t,\text{\ \ an impossibility.}$$ 

\noindent {\bf Case 2:  $r = 1$ or $r=2$.}
\vspace{.03in}
\newline  Our assumption is that $\fo(\ctr)\leq 2t+r-1$.  Therefore, we may choose a solvable distribution $D$ of $\ctr$ of size $2t+r-1$.  In each case, we modify $D$ to create a solvable distribution on a smaller cycle graph with fewer than the number of pebbles than we know is required.

{\it Subcase 2.1}: $\ctr$ contains a vertex $v$ with $D(v)=1$.
\vspace{.03in}
\newline The proof for  Subcase 2.1 of Theorem 1 adapts directly.

{\it Subcase 2.2}: $D$ places exactly 2 pebbles on each occupied vertex of $\ctr$.
\vspace{.03in}
\newline First,  since $|D|=2t+r-1$ and $|D|$ is even, we have $r=1$.  Number the vertices of $\co$ sequentially, $\co=(v_1, v_2, \dots, v_{3t+1})$, and consider the corresponding sequence of the number of pebbles on the vertices of $\co$.  Note that we may assume that there are at most two consecutive unoccupied vertices in $D$, since $D$ was assumed to be solvable.  Also, there must be a subsequence of vertices, $v_i, v_{i+1}, v_{i+2}$, with $D([v_i, v_{i+1}, v_{i+2}])=[2, 0, 2]$ or $D([v_i, v_{i+1}, v_{i+2}])=[2, 2, 0]$.  Otherwise, there would be {\it exactly} two unoccupied vertices between every pair of occupied vertices, yielding $r=0$, a contradiction.

In each case, to obtain a new distribution, $D^*$, remove vertices $v_{i+1}$ and $v_{i+2}$ and their associated pebbles.  In the first case,  $D^*$ is a solvable distribution since $v_{i+3}$ is either occupied or can be pebbled by $v_i$, and vertices $v_{i-1}$ and $v_{i+4}$ are unaffected.  In the second case, no vertices are affected by the removal of $v_{i+1}$ and $v_{i+2}$, so $D^*$ is solvable.  Also $D^*$ is a distribution on $\cm$ with $|D^*|=2t-2$.  We have reached a contradiction since by hypothesis,  $\fo(\cm)=\fo(C_{3(t-1)+2})=2(t-1)+2=2t$. 

{\it Subcase 2.3}: $\ctr$ contains some vertex $v_1$ such that $D(v_1)\geq 3$.
\vspace{.03in}
\newline We consider how to construct $D^*$ in each of three possible cases.  For the first two cases, we number the vertices of $\ctr$ sequentially, $\ctr=(v_1, v_2, \dots, v_{3t+r})$ in either direction.  If either of the vertices adjacent to $v_1$, say $v_2$, is unoccupied, $D^*$ may be constructed by removing $v_2$ and 2 pebbles from $v_1$, and then adding 1 pebble to $v_{3t+r}$.  $D^*$ permits at least as many pebbles to move to $v_3$ and $v_{3t+r}$ as does $D$, therefore $D^*$ is solvable.  Say both vertices adjacent to $v_1$ are occupied.  If $D(v_1)=3$, remove $v_1$ and its pebbles and place a pebble on each vertex adjacent to $v_1$ to obtain the distribution $D^*$.  $D^*$ permits at least as many pebbles to move to $v_2$ and $v_{3t+r}$ as does $D$, so $D^*$ is solvable.  

Finally, if $D(v_1)>3$, locate the unoccupied vertex closest to $v_1$ and number the vertices of $\ctr$ sequentially in that direction.  Designate the unoccupied vertex closest to $v_1$ to be $v_j$.  To obtain $D^*$, remove 3 pebbles from $v_1$, place 2 additional pebbles on $v_{3t+r}$ and remove vertex $v_j$.  Observe that $D^*$ permits as least as many pebbles to move to $v_{3t+r}$ as does $D$.  Finally, note that the vertices $v_2, v_3,\dots,v_{j-1}$ are occupied in both $D$ and $D^*$.    The 3 pebbles that were removed from $v_1$ could contribute at most 2 pebbles to vertex $v_{j-1}$ and therefore none to vertex $v_{j+1}$ under distribution $D$.  Thus, $D^*$ is indeed solvable.

Thus whenever $D$ places at least three pebbles on some vertex of $\ctr$, we can create a solvable distribution, $D^*$, on $C_{3t+r-1}$ with $|D^*|=2t+r-2$.  This contradicts our hypothesis.

Having reached a contradiction in all cases, we conclude that our initial assumption is false, and the theorem holds for all positive integers exceeding 2.

\section{Graham's conjecture}

An important open question in (non-optimal) graph pebbling is Graham's conjecture, a statement about the pebbling number of graph products.  We define the product of two graphs $G$ and $H$ to be the graph with vertex set $V(G) \times V(H)$ (Cartesian product) and with edge set
\begin{align}
E(G\times H) = &\{\bigl( (g_1,h_1),(g_2,h_2) \bigr) \vert\, g_1=g_2 \text{ and } (h_1,h_2)\in V(H)\} \notag \\
&\cup \{\bigl( (g_1,h_1),(g_2,h_2) \bigr) \vert\, h_1=h_2 \text{ and } (g_1,g_2)\in V(G)\} \notag  .
\end{align}
Recall that $f(G)$ denotes the pebbling number of the graph $G$.
\begin{conjecture} \cite{Chung} For any two graphs $G$ and $H$, $f(G\times H) \leq f(G)f(H)$.
\end{conjecture}
While this question is still open, many successful results in support have appeared.  (See, for example, \cite{Chung}, \cite{Moews92}, and \cite{SnevilyFoster}.)  Most of these take the form of showing that $f(G\times H) \leq f(G)f(H)$ for $G$ from a particular family of graphs and $H$ satisfying the 2-pebbling property\footnote{Given a distribution $D$ on a graph $G$ with $|D|=p$ and $D(v)>0$ for exactly $q$ vertices $v$, $G$ satisfies the 2-pebbling property if $p+q>2f(G)$.}.  (Herscovici introduces a similar property\footnote{$G$ satisfies the path property if for every path $P_m$, $f(P_m \times G)\leq f(P_m)f(G)$.}, as there are known to be graphs -- the Lemke graphs \cite{SnevilyFoster} -- that do not satisfy the 2-pebbling property \cite{Herscovici}.)

The analogous statement for optimal pebbling numbers is easier.  \begin{theorem}\label{Graham} For any two graphs $G$ and $H$, $\fo(G\times H) \leq \fo(G)\fo(H)$\footnote{Hung-Lin Fu and Chin-Lin Shiue state this result in \cite{FuShiue}.  The proof has yet to appear;  the status of the paper containing their proof is unclear \cite{Fu}.}.
\end{theorem}
{\bf Proof.}
Let $D_G$ and $D_H$ be optimally-sized solvable distributions on $G$ and $H$, respectively.  Define the distribution $D$ on $G \times H$ by placing $D_G(v)D_H(w)$ pebbles on each vertex $(v,w) \in V(G \times H)$.  Note that $|D|= \fo(G)\fo(H)$, so the theorem is proven once we establish $D$'s solvability.

An arbitrary vertex $(v,w) \in V(G \times H)$ lies in the subgraph $\{v\}\times H$, which is isomorphic to $H$.  Likewise, each vertex $(v,w^\prime)$ of $\{v\}\times H$ lies in a subgraph $G \times \{w^\prime\}$ isomorphic to $G$.  Now, consider the distribution $D_G$ on $G$.  Suppose $j$ pebbles can be moved to $v$ under $D_G$.  Note that $j\geq1$ since $D_G$ is a solvable distribution.  Then, starting from the distribution $D$ on $G\times H$, we may move at least $jD_H(w^\prime)$ pebbles to the vertex $(v, w^\prime)$.  In particular, we may move at least $D_H(w^\prime)$ pebbles to each vertex $(v,w^\prime) \in V(\{v\}\times H)$.  But, since $D_H$ is solvable on $H$ and we have at least $|D_H|$ pebbles arranged appropriately on the vertices of $\{v\}\times H$, we can use edges of the form $\bigl((v,w^\prime),(v,w^{\prime \prime})\bigr)$ to move a pebble to $(v,w)$.  Since the vertex $(v, w)$ was arbitrary, the distribution $D$ is solvable, thus proving the theorem.


Equality in Theorem \ref{Graham} is achieved by $P_3 \times P_3$, which has optimal pebbling number 4. (Label the vertices of $P_3$ sequentially as $v_1, v_2, v_3$ and place all 4 pebbles on $(v_2,v_2) \in V(P_3 \times P_3)$.)    
The inequality, however, is required: observe that $C_4$, with optimal pebbling number 3, is isomorphic to $P_2 \times P_2$, and $\fo(P_2)\fo(P_2)=2\cdot 2 = 4$.

Establishing an upper bound for the optimal pebbling number of a graph requires one only to prove that an appropriately-sized distribution is solvable.  Similarly, lower bounds for (non-optimal) pebbling numbers can be proven by demonstrating the non-solvability of a distribution of size one less than the claimed lower bound.  It is the constructive nature of demonstrating upper bounds for optimal pebbling numbers that leads to the ease of proving the optimal-pebbling version of Graham's conjecture.  This suggests that investigating lower bounds for optimal pebbling numbers may be of interest in future research.

\section{Acknowledgment}
This work was carried out while the first author was a guest of the Department of Computer Science, Mathematics, and Statistics at Mesa State College.  The hospitality of the department is remembered with gratitude!

\end{document}